\documentclass{amsart}

\usepackage[T1]{fontenc}
\usepackage[utf8]{inputenc}
\usepackage{lmodern}
\usepackage{amsmath,amssymb,mathtools}
\usepackage{booktabs,array,longtable}
\usepackage{microtype}
\usepackage{xurl}
\usepackage[hidelinks]{hyperref}

\allowdisplaybreaks[2]
\providecommand{\tightlist}{%
  \setlength{\itemsep}{0pt}\setlength{\parskip}{0pt}}

\hypersetup{
  pdftitle={An independent proof of the even-dimensional S-matrix inequality},
  pdfauthor={Yinjie Li}}

\title{An independent proof of the even-dimensional S-matrix inequality}
\DeclareRobustCommand{\AuthorAffiliationMarker}{a}
\DeclareRobustCommand{\AuthorAffiliation}{\normalfont\textsuperscript{a} Independent Researcher, Nanjing, China}
\DeclareRobustCommand{\AuthorEmail}{\normalfont Email: \texttt{stewarteur@gmail.com}}
\author[Yinjie Li]{Yinjie Li\textsuperscript{\AuthorAffiliationMarker}\\
\AuthorAffiliation\\
\AuthorEmail}

\begin{document}
\begin{abstract}
Harwit and Sloane conjectured that every nonsingular entrywise-nonnegative
matrix $A\in\mathbb R^{n\times n}$ satisfies
$\|A^{-1}\|_F\ge 2n(n+1)^{-1}\|A\|_{\max}^{-1}$, with equality
precisely for positive multiples of $S$-matrices. Zhang has given a complete
proof of this conjecture by a centered pseudoinverse and spectral variance
method. We present an independently obtained, structurally different proof
of the strict even-dimensional inequality. Starting from the structural
identities of Frankel and Urschel, we derive an exact global defect budget
and combine binary rounding, fixed intersections, and Gram projection.
A ten-row obstruction handles every even $n\ge66$; a finite exact
calculation handles $4\le n\le64$, $n\ne6$; and a multi-column energy
argument treats $n=6$. The order-two case is elementary. The even-dimensional
argument is formalized in Lean 4, conditional on Frankel--Urschel Lemma 2.1
as an explicit external mathematical input. The finite evaluations use
Lean's native evaluator; their trust boundary and exact certificates are
documented. Together with Cheng's odd-dimensional theorem, the argument
recovers the full S-matrix theorem and its equality characterization.
\end{abstract}
\maketitle

\hypertarget{introduction-and-main-theorem}{%
\section{Introduction and main
theorem}\label{introduction-and-main-theorem}}

For a real matrix \(A=(A_{ij})\), write

\[
\|A\|_F=\sqrt{\operatorname{trace}(A^TA)},
\qquad
\|A\|_{\max}=\max_{i,j}|A_{ij}|,
\]

and let \(\mathbf 1\) denote the all-ones column vector. An \(S\)-matrix
of order \(n\) is a matrix \(S\in\{0,1\}^{n\times n}\) satisfying

\[
S^TS=\frac{n+1}{4}(I+\mathbf 1\mathbf 1^T).
\]

The problem arose in Hadamard transform optics. Harwit and Sloane asked
for the least possible Frobenius norm of the inverse of a nonsingular
matrix with entries in \([0,1]\) and conjectured that the optimum is
attained exactly by \(S\)-matrices {[}4,5{]}. Cheng proved the
conjecture for odd \(n\) and obtained the weaker strict estimate

\[
\|A^{-1}\|_F>
\frac{2\sqrt{n^2-2n+2}}{n}\|A\|_{\max}^{-1}
\]

for even \(n\ge4\) {[}1{]}. Drnovšek later gave a shorter proof of these
bounds {[}2{]}. Frankel and Urschel introduced a stability analysis of
the even case and proved the conjectured estimate for all \(n\ge1000\)
{[}3{]}, while noting that their numerical threshold was not optimized.

Zhang has given a complete proof of the conjecture in all dimensions
{[}6{]}. His even-dimensional argument works with a centered matrix whose
rows lie in a projected cube: parity estimates for the sign Gram matrix,
a variational spectral variance bound, and an exact inverse decomposition
give a centered pseudoinverse estimate and recover the mean direction.

The present proof was obtained independently of Zhang's argument. It
develops the Frankel--Urschel structural identities into an exact defect
budget, then uses quantitative near-integrality, binary rounding, fixed row
weights and intersections, and combinatorial Gram rigidity. Thus
Frankel--Urschel provide the structural starting point and a large-dimension
theorem; our method completes the even-dimensional cases by discrete
rigidity and exact energy estimates; Zhang supplies a complete proof through
centered spectral analysis. The contribution of this paper is the independent
even-dimensional proof and its formal verification, rather than a claim to
the first complete resolution of the conjecture. The theorem below records
the common conclusion; only its even-dimensional part is proved here by the
new argument, while the odd-dimensional part and equality classification
are imported from Cheng's theorem.

\textbf{Theorem 1.1 (S-matrix conjecture).} Let
\(A\in\mathbb R^{n\times n}\) be nonsingular and entrywise nonnegative.
Then

\[
\boxed{
\|A^{-1}\|_F\ge \frac{2n}{n+1}\,\|A\|_{\max}^{-1}
}.
\]

Equality holds if and only if \(A\) is a positive multiple of an
\(S\)-matrix.

There is no \(S\)-matrix of even order, since such a matrix would
produce a Hadamard matrix of odd order \(n+1>1\). Cheng's theorem
settles every odd order, so it remains to establish the strict estimate
for even \(n\). The order \(n=2\) is elementary and is treated in
Section 8. Our even-dimensional argument for \(n\ge4\) has three parts.

\begin{enumerate}
\def\labelenumi{\arabic{enumi}.}
\tightlist
\item
  For \(n\ge66\), ten low-defect rows round to a binary matrix whose
  Gram matrix has two eigenvalues. Minimizing the induced quadratic form
  over a box then contradicts the available column energy.
\item
  For \(4\le n\le64\), \(n\ne6\), we first determine the rounded row
  weights and pairwise intersections. The remaining Gram inequalities
  form a finite family that is checked by exact rational arithmetic.
\item
  At \(n=6\), the estimates from a single column are insufficient. The
  global defect budget nevertheless forces five rounded triples, whose
  exceptional intersections create more energy across several columns of
  the Frankel--Urschel error matrix than the budget permits.
\end{enumerate}

No floating-point estimate enters the proof. The finite calculation and
its connection to the matrix argument are formalized in Lean; Section 7
briefly describes the scope of this verification. Notation specific to
the proof is introduced where it is first used.

\hypertarget{the-frankelurschel-structural-lemma}{%
\section{The Frankel--Urschel structural
lemma}\label{the-frankelurschel-structural-lemma}}

The following result is the sole external mathematical input to the new
even-dimensional argument; it is Lemma 2.1 of Frankel and Urschel
{[}3{]}. Let \(n>2\) be even, and suppose that
\(B\in\mathbb R^{n\times n}\) is nonsingular, entrywise nonnegative, and
satisfies

\[
\|B\|_{\max}\le1,
\qquad
\|B^{-1}\|_F\le\frac{2n}{n+1}.
\]

Put \(s=1/[2(n-2)]\) and define

\[
r=B\mathbf1-\frac{(n-1)^2}{2(n-2)}\mathbf1
=B\mathbf1-\left(\frac n2+s\right)\mathbf1,
\tag{2.1}
\]

and

\[
\begin{aligned}
H(B)=\frac{\sqrt{n(n-2)}}{2(n-1)}
\Bigg[&
\sqrt{\frac n{n-2}}\frac n2 B^{-1}\\
&-\sqrt{\frac{n-2}{n}}
\left(
\frac{2(n-1)}{n-2}I-\frac2n\mathbf1\mathbf1^T
\right)B^T
\Bigg].
\end{aligned}
\tag{2.2}
\]

Set

\[
\begin{aligned}
h(B)&=
\left\|
\sqrt{\frac{n-2}{n}}
\left(
\frac{2(n-1)}{n-2}I-\frac2n\mathbf1\mathbf1^T
\right)B^T
\right\|_F^2,\\
c&=\frac{n(n^2-2n+2)}{n-2}-h(B).
\end{aligned}
\]

The lemma states that \(0\le c<1\) and gives the following three
conclusions:

\[
\|r\|_2^2+
\frac{(n-1)^2}{n-2}
\sum_{i,j=1}^nB_{ij}(1-B_{ij})
=\frac{cn}{4}+\frac{n}{4(n-2)^2},
\tag{2.3}
\]

\[
\|H(B)\|_F^2
\le
\frac{n(n-2)}{4(n-1)^2}
\left[
\frac{n(n^2-2n-2)}{(n-2)(n+1)^2}-c
\right],
\tag{2.4}
\]

and the exact matrix identity

\[
\begin{aligned}
B(B^T+H(B))={}&
\frac{n^2}{4(n-1)}I
+\frac{(n-1)^3}{4n(n-2)}\mathbf1\mathbf1^T\\
&+\frac{n-1}{2n}
(r\mathbf1^T+\mathbf1r^T)
+\frac{n-2}{n(n-1)}rr^T.
\end{aligned}
\tag{2.5}
\]

Equations (2.3), (2.4), and (2.5) are respectively conclusions (1), (2),
and (3) of {[}3, Lemma 2.1{]}.

For completeness, (2.5) follows directly from (2.2). If

\[
M=\frac{2(n-1)}{n-2}I-\frac2n\mathbf1\mathbf1^T,
\]

then

\[
\begin{aligned}
B(B^T+H(B))
&=BB^T+\frac{n^2}{4(n-1)}I
-\frac{n-2}{2(n-1)}BMB^T\\
&=\frac{n^2}{4(n-1)}I
+\frac{n-2}{n(n-1)}(B\mathbf1)(B\mathbf1)^T.
\end{aligned}
\]

Substituting \(B\mathbf1=r+(n/2+s)\mathbf1\) gives (2.5), including the
positive coefficient of \(rr^T\). The negative sign in one intermediate
display of {[}3, Lemma 2.1(3){]} is therefore typographical; the
statement of that lemma has the positive sign. We shall abbreviate

\[
\beta=\frac{(n-1)^3}{4n(n-2)}.
\]

For a row \(i\), put

\[
u_i=\sum_{j=1}^n B_{ij}(1-B_{ij}),
\qquad
h_i^2=\sum_{j=1}^n H(B)_{ji}^2,
\qquad
\gamma=\frac{n-2}{(n-1)^2}.
\]

For later use, set

\[
K=\frac{4(n-1)^2}{n(n-2)},
\qquad
\mathcal D_i=\frac4n r_i^2+K(u_i+h_i^2),
\qquad
\mathcal U_i=\frac{\mathcal D_i}{K}
=u_i+h_i^2+\gamma r_i^2.
\]

Multiplying (2.3) by \(4/n\) and (2.4) by \(K\), then adding the two
inequalities, eliminates \(c\). Hence

\[
\sum_{i=1}^n\mathcal D_i\le\tau_n,
\tag{2.6}
\]

where

\[
\boxed{
\tau_n=\frac1{(n-2)^2}
+\frac{n(n^2-2n-2)}{(n-2)(n+1)^2}
}.
\]

For \(n\ge6\), this number is less than one: indeed,

\[
1-\tau_n=
\frac{2n^3-6n^2-2n+3}{(n-2)^2(n+1)^2}>0.
\]

The defect bound (2.6) and the matrix identity (2.5) are the common
starting point for the three even-dimensional arguments.

\hypertarget{normalization-and-contradiction-setup}{%
\section{Normalization and contradiction
setup}\label{normalization-and-contradiction-setup}}

Let \(n\ge4\) be even and suppose, for contradiction, that \(A\) is a
nonsingular entrywise nonnegative matrix for which

\[
\|A^{-1}\|_F\le\frac{2n}{n+1}\|A\|_{\max}^{-1}.
\]

Put

\[
M=\|A\|_{\max}>0,
\qquad
B=M^{-1}A.
\]

Then \(0\le B_{ij}\le1\), \(\|B\|_{\max}=1\), and \(B^{-1}=MA^{-1}\).
Thus

\[
\|B^{-1}\|_F\le\frac{2n}{n+1},
\tag{3.1}
\]

so the conclusions of Section 2 apply.

Round every entry of \(B\) to a nearest element of \(\{0,1\}\); the
choice at a tie is immaterial. Write

\[
\begin{gathered}
C=\operatorname{round}(B)\in\{0,1\}^{n\times n},
\qquad E=B-C,\\
\epsilon_i=\sum_jE_{ij},
\qquad
\sigma_i=\sum_j|E_{ij}|,
\qquad
w_i=\sum_jC_{ij}.
\end{gathered}
\]

For every entry, let

\[
\eta_{ij}=|E_{ij}|\le\frac12.
\]

Since \(B_{ij}(1-B_{ij})=\eta_{ij}(1-\eta_{ij})\),

\[
u_i=\sum_j\eta_{ij}(1-\eta_{ij}),
\qquad
\sigma_i\le2u_i,
\qquad
u_i\ge\sigma_i(1-\sigma_i).
\tag{3.2}
\]

The last inequality follows from

\[
u_i=\sigma_i-\sum_j\eta_{ij}^2
\ge\sigma_i-\sigma_i^2.
\]

Combining (2.1) with \(B=C+E\) gives the row-sum relation

\[
\boxed{
r_i=\left(w_i-\frac n2\right)+\epsilon_i-s.
}
\tag{3.3}
\]

Thus \(r_i=-s+\epsilon_i\) once the rounded row has weight \(n/2\).

\hypertarget{large-even-dimensions-nge66}{%
\section{\texorpdfstring{Large even dimensions:
\(n\ge66\)}{Large even dimensions: n\textbackslash ge66}}\label{large-even-dimensions-nge66}}

In this section fix

\[
t=9,
\qquad
m=t+1=10.
\]

This choice retains enough rows to force a useful Gram projection while
leaving sufficient defect control at the limiting order \(n=66\).

The estimates below require \(t\ge4\) and \(n\ge4t\), both of which hold
for every even \(n\ge66\).

\hypertarget{low-defect-rows-and-their-gram-matrix}{%
\subsection{Low-defect rows and their Gram
matrix}\label{low-defect-rows-and-their-gram-matrix}}

We first select ten distinct indices

\[
i_1,\ldots,i_{10}
\]

such that, for every selected index \(i\),

\begin{align}
w_i=[C\mathbf1]_i&=\frac n2, \tag{4.1}\\
h_i^2&<\frac1{4(n-t)}, \tag{4.2}\\
\sigma_i=\sum_{j=1}^n|B_{ij}-C_{ij}|&<e,
\qquad e:=\frac1{4(n-t)-2}, \tag{4.3}\\
|r_i|&<\frac1{2(n-2)}+e. \tag{4.4}
\end{align}

Indeed, (2.3)--(2.4) imply

\[
\frac4n\|r\|_2^2
+K\sum_i u_i
+K\|H(B)\|_F^2
\le
\frac1{(n-2)^2}
+\frac{n(n^2-2n-2)}{(n-2)(n+1)^2}<1
\tag{4.5}
\]

for \(n\ge6\), by the formula following (2.6). If fewer than \(t+1\)
indices satisfied \(\mathcal D_i<1/(n-t)\), then at least \(n-t\) terms
in (4.5) would sum to at least \(1\). Hence the selected indices may be
chosen so that \(\mathcal D_i<1/(n-t)\). Since \(K>4\), this gives (4.2)
and \(u_i<1/[4(n-t)]\). In particular, each \(\eta_{ij}<1/[2(n-t)]\),
because \(\eta_{ij}(1-\eta_{ij})\ge \eta_{ij}/2\). Therefore

\[
\sum_j|B_{ij}-C_{ij}|
\le
\frac{1}{1-1/[2(n-t)]}\sum_jB_{ij}(1-B_{ij})
<\frac1{4(n-t)-2}.
\]

Moreover,

\[
|[B\mathbf1]_i-\operatorname{round}([B\mathbf1]_i)|
\le\sigma_i<\frac1{4(n-t)-2}.
\]

The same defect bound gives

\[
|r_i|<\frac12\sqrt{\frac n{n-t}}.
\]

For \(t=9\) and \(n\ge66\), we have

\[
1-s-e>\frac12\sqrt{\frac n{n-t}}.
\]

For example, \(s\le1/128\), \(e\le1/226\), and
\(\sqrt{n/(n-9)}\le\sqrt{66/57}<11/10\), so the left-hand side exceeds
\(9/10\) while the right-hand side is below \(11/20\).

If \(w_i\ne n/2\), then (3.3) and \(|\epsilon_i|\le\sigma_i<e\) would
instead give \(|r_i|>1-s-e\), a contradiction. Thus (4.1) holds, and
(3.3) immediately yields (4.4).

Let \(\widehat C\in\{0,1\}^{10\times n}\) be the restriction of \(C\) to
these rows. Put

\[
q=\left\lceil\frac n4\right\rceil,
\qquad
\rho=\left\lfloor\frac n4\right\rfloor.
\]

For distinct selected rows, expand \(B=C+E\) in (2.5). The terms
involving \(E\) are bounded by (4.3), the row-sum terms by (4.4), and
\(|B_iH(B)_{*,j}|\) by Cauchy--Schwarz together with (4.2). This gives

\[
\begin{aligned}
\left|
(CC^T)_{ij}-\frac{(n-1)^3}{4n(n-2)}
\right|
<&\frac1{2(n-t)-1}
+\frac1{2(n-2)}\\
&+\frac1{4(n-t)-2}
+\frac12\sqrt{\frac n{n-t}}.
\end{aligned}
\]

For \(t=9\) and \(n\ge66\), the right-hand side is less than \(7/10\).
Indeed, the four terms decrease with \(n\), and at \(n=66\) the bound
\(\sqrt{66/57}<11/10\) gives

\[
\frac1{113}+\frac1{128}+\frac1{226}+\frac{11}{20}<\frac35<\frac7{10}.
\]

Also

\[
\left|
\frac{(n-1)^3}{4n(n-2)}
-\operatorname{round}\left(
\frac{(n-1)^3}{4n(n-2)}
\right)
\right|<\frac3{10},
\]

and the indicated nearest integer is \(\rho\). Since \((CC^T)_{ij}\) is
an integer, the two strict inequalities imply

\[
(\widehat C\widehat C^T)_{ij}=\rho
\qquad(i\ne j).
\tag{4.6}
\]

Together with the row weight \(n/2=q+\rho\), this gives

\[
\boxed{
\widehat C\widehat C^T=qI_{10}+\rho J_{10}.
}
\]

\hypertarget{perturbation-and-gram-projection}{%
\subsection{Perturbation and Gram
projection}\label{perturbation-and-gram-projection}}

Define

\[
\boxed{
\mathcal E(n,t)=
\frac1{2(n-t)-1}
+\frac1{2(n-2)}
+\frac1{4(n-t)-2}
+\frac1{(8(n-t)-4)\sqrt{n-t}}.
}
\tag{4.7}
\]

Comparing \(C(C^T+H(B))\) with the right-hand side of (2.5), and
bounding the centered correction and the three rounding terms, gives for selected
\(i\ne j\)

\[
\begin{aligned}
&\left|
[C(C^T+H(B))]_{ij}
-\frac{(n-1)^3}{4n(n-2)}
\right|\\
&\qquad<
2e+\left(e+\frac1{2(n-2)}\right)
+\frac{e}{2\sqrt{n-t}}
=\mathcal E(n,t),
\end{aligned}
\tag{4.8}
\]

Here is the error allocation behind (4.8). Put \(a_1=(n-1)/(2n)\),
\(a_2=(n-2)/(n(n-1))\), \(h=1/(2\sqrt{n-t})\), and \(R=s+e<1\).
The exact off-diagonal error is

\[
\begin{aligned}
&[C(C^T+H(B))]_{ij}-\beta\\
&\qquad=a_1(r_i+r_j)+a_2r_ir_j
-E_iC_j^T-B_iE_j^T-E_iH(B)_{*,j}.
\end{aligned}
\]

Since \(|r_i|,|r_j|<R\), \(a_1,a_2>0\), and

\[
2a_1+a_2=1-\frac1{n(n-1)}<1,
\]

we have

\[
|a_1(r_i+r_j)+a_2r_ir_j|
<2a_1R+a_2R^2<R=s+e.
\]

The remaining three terms have absolute values less than \(e,e,eh\),
respectively. Thus the full error is less than
\(3e+s+eh=\mathcal E(n,t)\), proving (4.8).
On the diagonal, the centered correction \(2a_1r_i+a_2r_i^2\) is bounded
by \(|r_i|<e+s\); the remaining three terms have bounds \(e,e,eh\).
Consequently, on the diagonal

\[
\begin{aligned}
&\left|
[C(C^T+H(B))]_{ii}
-\left(
\frac{n^2}{4(n-1)}
+\frac{(n-1)^3}{4n(n-2)}
\right)
\right|\\
&\qquad<\mathcal E(n,t).
\end{aligned}
\tag{4.9}
\]

Here

\[
2e=\frac1{2(n-t)-1},
\qquad
\frac{e}{2\sqrt{n-t}}
=\frac1{(8(n-t)-4)\sqrt{n-t}}.
\]

Let

\[
y=(-1)^{n/2+1}H(B)_{*,i_{10}}\in\mathbb R^n,
\qquad
z=\widehat C y\in\mathbb R^{10}.
\tag{4.10}
\]

By (4.2),

\[
\boxed{
\|y\|_2^2<\frac1{4(n-t)}.
}
\tag{4.11}
\]

Apply (4.8)--(4.9) to the column indexed by \(i_{10}\), subtract the
Gram entries from (4.6), and use the sign in (4.10). We obtain

\begin{align}
|z_j-\tfrac14|&<\delta_1 \qquad(j=1,\ldots,t), \notag\\
|z_{t+1}|&<\delta_0, \notag\\
\delta_1&=\frac{n-1}{4(n-2)n}+\mathcal E(n,t), \tag{4.12}\\
\delta_0&=\frac{2n^2-4n+1}{4(n-2)(n-1)n}+\mathcal E(n,t). \tag{4.13}
\end{align}

The rational terms in (4.12)--(4.13) are the distances between the
centers in (4.8)--(4.9) and the corresponding integer Gram entries.

We use the following elementary projection lemma.

\textbf{Lemma 4.1.} Let \(m=t+1\), let \(D\in\mathbb R^{m\times n}\)
satisfy

\[
DD^T=qI_m+\rho J_m,
\qquad q>0,
\quad \rho\ge0,
\]

and suppose \(z=Dy\) satisfies

\[
z_1,\ldots,z_t\in[a,a+2\delta_1],
\qquad
z_m\in[-d,d],
\tag{4.14}
\]

where

\[
a=\frac14-\delta_1,
\qquad
d=\delta_0.
\]

Assume

\[
a(q+\rho)>\rho d,
\qquad
d(q+\rho t)<\rho ta.
\tag{4.15}
\]

Then

\[
\boxed{
\|y\|_2^2\ge
\frac{t(a-d)^2}{(t+1)q}
+\frac{(ta+d)^2}{(t+1)(q+\rho(t+1))}.
}
\tag{4.16}
\]

\textbf{Proof.} Since \(D\) has full row rank, the least-norm vector
satisfying \(Dy=z\) is

\[
y_*=D^T(DD^T)^{-1}z.
\]

Hence

\[
\|y\|_2^2\ge z^T(DD^T)^{-1}z.
\tag{4.17}
\]

The inverse Gram matrix is

\[
(qI_m+\rho J_m)^{-1}
=\frac1qI_m-
\frac{\rho}{q(q+\rho m)}J_m.
\]

Thus the right-hand side of (4.17) is the convex quadratic form

\[
\Phi(z)=
\frac1q\sum_{i=1}^m z_i^2
-\frac{\rho}{q(q+\rho m)}
\left(\sum_{i=1}^m z_i\right)^2.
\tag{4.18}
\]

At \(z_*=(a,\ldots,a,d)\), the derivatives in the first \(t\)
coordinates have the sign of \(a(q+\rho)-\rho d\), whereas the
derivative in the last coordinate has the sign of
\(d(q+\rho t)-\rho ta\). Conditions (4.15) therefore make the first
\(t\) coordinates minimal at their lower endpoints and the last
coordinate minimal at its upper endpoint. Convexity now shows that
\(z_*\) minimizes \(\Phi\) throughout the box (4.14). Substituting
\(z_*\) into (4.18) and separating the all-ones direction from its
orthogonal complement gives

\[
\Phi(z_*)=
\frac{t(a-d)^2}{(t+1)q}
+\frac{(ta+d)^2}{(t+1)(q+\rho(t+1))}.
\]

This proves the lemma. \(\square\)

For \(n\ge66\) and \(t=9\), the bounds

\[
a>0.22,
\qquad d<0.03,
\qquad \frac q\rho\le\frac{17}{16}
\]

imply both inequalities in (4.15). With

\[
Q(n)=\frac{9(a-d)^2}{10q}+\frac{(9a+d)^2}{10(q+10\rho)},
\]

(4.11) and (4.16) are incompatible as soon as

\[
Q(n)>\frac1{4(n-9)}.
\tag{4.19}
\]

\hypertarget{endpoint-estimates}{%
\subsection{Endpoint estimates}\label{endpoint-estimates}}

We first treat the critical order \(n=66\).

For \(n=66\),

\[
q=17,
\qquad
\rho=16,
\qquad
n-t=57.
\]

Using \(\sqrt{57}>15/2\),

\[
\mathcal E(66,9)<\overline{\mathcal E}_{66}
=\frac1{113}+\frac1{128}+\frac1{226}+\frac1{3390}
=\frac{4639}{216960}.
\tag{4.20}
\]

Because the estimate in (4.20) is strict, substituting
\(\overline{\mathcal E}_{66}\) into (4.12)--(4.13) gives

\[
a>\bar a:=\frac{2145719}{9546240},
\qquad
d<\bar d:=\frac{721649}{24820224}.
\]

Within the region (4.15), the lower bound (4.16) is increasing in \(a\)
and decreasing in \(d\). Therefore

\[
Q(66)>
\frac9{170}(\bar a-\bar d)^2
+\frac1{1770}(9\bar a+\bar d)^2.
\]

Exact rational arithmetic gives

\[
\begin{aligned}
&\frac9{170}(\bar a-\bar d)^2
+\frac1{1770}(9\bar a+\bar d)^2
-\frac1{228}\\
&\qquad=
\frac{9018303786068791}
{440247800602489651200}>0.
\end{aligned}
\]

Since \(1/228=1/[4(66-9)]\), this proves (4.19) at \(n=66\),
contradicting (4.11).

For the remaining even orders \(n\ge68\), put

\[
b_1(n)=\frac{n-1}{4(n-2)n},
\qquad
b_0(n)=\frac{2n^2-4n+1}{4(n-2)(n-1)n}.
\]

Then

\[
a-d=\frac14-b_1(n)-b_0(n)-2\mathcal E(n,9),
\tag{4.21}
\]

\[
9a+d=\frac94-9b_1(n)+b_0(n)-8\mathcal E(n,9).
\tag{4.22}
\]

The functions \(b_1(n)\), \(b_0(n)\), and \(\mathcal E(n,9)\) are
decreasing for \(n\ge68\). In particular,

\[
b_1'(n)=-\frac{(n-1)^2+1}{4n^2(n-2)^2}<0,
\]

and

\[
b_0'(n)=-
\frac{2(n-2)^4+8(n-2)^3+11(n-2)^2+6(n-2)+2}
{4n^2(n-1)^2(n-2)^2}<0.
\]

Every summand of \(\mathcal E(n,9)\) is also decreasing.

Using \(\sqrt{59}>192/25\), define

\[
\overline{\mathcal E}_{68}
=\frac1{117}+\frac1{132}+\frac1{234}
+\frac{25}{468\cdot192}
=\frac{20435}{988416}.
\tag{4.23}
\]

Equations (4.21)--(4.23) imply, for every \(n\ge68\),

\[
a-d>A_0,
\qquad
A_0=\frac{111148175}{562902912},
\]

and, after discarding the positive term \(b_0(n)\) in (4.22),

\[
9a+d>B_0,
\qquad
B_0=\frac{2153959}{1050192}.
\]

For even \(n\),

\[
q\le\frac{n+2}{4},
\qquad
\rho\le\frac n4,
\qquad
q+10\rho\le\frac{11n+2}{4}.
\]

Therefore

\[
Q(n)>
\frac{36A_0^2}{10(n+2)}
+\frac{4B_0^2}{10(11n+2)}.
\]

After multiplication by \(4(n-9)\), the right-hand side is a positive
linear combination of

\[
\frac{n-9}{n+2}
\qquad\text{and}\qquad
\frac{n-9}{11n+2},
\]

both strictly increasing for \(n>9\). It is therefore enough to check
\(n=68\). Exact rational arithmetic gives

\[
\begin{aligned}
&4(68-9)
\left[
\frac{36A_0^2}{10(68+2)}
+\frac{4B_0^2}{10(11\cdot68+2)}
\right]-1\\
&\qquad=
\frac{2795300344016846513}
{1039695852359324160000}>0.
\end{aligned}
\]

Hence (4.19) holds for every even \(n\ge68\), again contradicting
(4.11).

The sharper error term in (4.7) matters only near the endpoint: with the
coarser estimate of {[}3{]}, the left-hand side of (4.19) at \(n=66\)
reaches only about \(0.99321\) of the required value.

\hypertarget{finite-even-dimensions-4le-nle64-nne6}{%
\section{\texorpdfstring{Finite even dimensions: \(4\le n\le64\),
\(n\ne6\)}{Finite even dimensions: 4\textbackslash le n\textbackslash le64, n\textbackslash ne6}}\label{finite-even-dimensions-4le-nle64-nne6}}

Fix an even \(n\) in this range. Here we use the full defect budget
(2.6), rather than replacing its right-hand side by \(1\). Retain the
notation \(q=\lceil n/4\rceil\) and \(\rho=\lfloor n/4\rfloor\).

Choose the number \(m=t+1\) of selected rows, and put \(d=n-t\), as
follows:

\[
\begin{array}{c|ccccc|c}
n&4&8&10&12&14&16,18,\ldots,64\\ \hline
t&1&3&5&5&7&n/2-2
\end{array}
\tag{5.1}
\]

The order \(n=6\) is excluded here and is treated independently in
Section 6. In particular, every parameter pair in (5.1) satisfies
\(d\ge n/2\), as required by Lemma 5.1.

The values in (5.1) were selected to maximize the exact certified margin
subject to this constraint; no optimality is claimed. There are thirty
parameter pairs: the thirty-one even orders from \(4\) through \(64\),
with \(n=6\) removed.

\hypertarget{rounding-to-a-fixed-gram-matrix}{%
\subsection{Rounding to a fixed Gram
matrix}\label{rounding-to-a-fixed-gram-matrix}}

Choose an index \(\ell\) and, for an integer \(t<n\), \(t\) further
indices such that

\begin{align}
\mathcal D_\ell&\le\frac{\tau_n}{n}, \notag\\
\mathcal D_i&\le\frac{\tau_n}{n-t}. \tag{5.2}
\end{align}

The first choice follows by averaging. For the second, if fewer than
\(t\) further indices satisfied (5.2), then at least \(n-t\) terms in
(2.6) would sum to at least \(\tau_n\). The distinguished index \(\ell\)
will supply the column of \(H(B)\). Put

\[
d=n-t,
\qquad
U=\frac{n(n-2)}{4(n-1)^2}\frac{\tau_n}{d}.
\tag{5.3}
\]

Then \(\mathcal U_i\le U\) for every selected row, while

\[
\mathcal U_\ell\le\frac{n(n-2)}{4(n-1)^2}\frac{\tau_n}{n}.
\tag{5.4}
\]

\textbf{Lemma 5.1.} Suppose \(d=n-t\ge n/2\), and a selected row
satisfies \(\mathcal U_i\le U\), with \(U\) as in (5.3). Then

\[
w_i=\frac n2.
\]

\textbf{Proof.} Put

\[
G_0=\gamma(1-s)^2.
\]

Using \(d\ge n/2\), we obtain

\[
G_0-U\ge
\frac{2n^4-4n^3-17n^2+18n+23}
{4(n-2)(n-1)^2(n+1)^2}>0.
\tag{5.5}
\]

The numerator is positive at \(n=4\); for \(n\ge5\) it equals

\[
n^2(2n^2-4n-17)+18n+23>0.
\]

Also \(G_0\le1/8\) for even \(n\ge4\), because after clearing
denominators this is equivalent to

\[
(n-4)(n^2-8n+13)\ge0.
\tag{5.6}
\]

Thus

\[
\sigma_i\le2u_i\le2U<\frac12.
\tag{5.7}
\]

If \(k=w_i-n/2\) is a nonzero integer, then \(|k-s|\ge1-s\). From (3.2),
(3.3), and \(\mathcal U_i\le U\),

\[
U\ge
\sigma_i(1-\sigma_i)
+\gamma(1-s-\sigma_i)^2
=:G(\sigma_i).
\]

The quadratic \(G\) is concave. Its minimum on \([0,1/2]\) is attained
at an endpoint, and

\[
G(0)=G_0,
\qquad
G(1/2)=\frac14+\gamma(1/2-s)^2>G_0.
\]

Hence \(U\ge G_0\), contradicting (5.5). Therefore \(w_i=n/2\).
\(\square\)

Consequently, (3.3) becomes

\[
\boxed{r_i=-s+\epsilon_i.}
\tag{5.8}
\]

\textbf{Lemma 5.2.} For every parameter pair \((n,t)\) in (5.1), any two
distinct selected rows have intersection \(\rho=\lfloor n/4\rfloor\);
that is,

\[
\langle C_i,C_j\rangle=\rho\qquad(i\ne j).
\]

\textbf{Proof.} For a fixed-weight selected row, first observe that

\[
U<s(1-s).
\tag{5.9}
\]

Using \(d\ge n/2\), the numerator of \(s(1-s)-U\), after clearing
positive denominators, is bounded below by

\[
7n^4-26n^3+10n^2+24n-1
=n^3(7n-26)+10n^2+24n-1>0.
\]

Equations (3.2), (5.7), and (5.9) imply \(\sigma_i<s\): otherwise the
monotonicity of \(x(1-x)\) on \([0,1/2]\) would give
\(u_i\ge s(1-s)>U\). Define

\[
\lambda=1-2\gamma s-2(1-\gamma)U,
\qquad
\delta=\frac{U-\gamma s^2}{\lambda}.
\tag{5.10}
\]

The preceding bounds also give \(\lambda>0\). Moreover,
\(|r_i|=|s-\epsilon_i|\ge s-\sigma_i\), so (3.2) yields

\[
\sigma_i(1-\sigma_i)+\gamma(s-\sigma_i)^2\le U.
\]

After expansion, and using \(\sigma_i\le2U\) in the quadratic term, this
gives

\[
U\ge\gamma s^2+\lambda\sigma_i,
\]

and therefore

\[
\sigma_i\le\delta.
\]

Set

\[
a_1=\frac{n-1}{2n},
\qquad a_2=\frac{n-2}{n(n-1)}.
\]

For distinct selected indices \(i,j\), (2.5) gives

\[
B_iB_j^T+B_iH(B)_{*,j}
=\beta+a_1(r_i+r_j)+a_2r_ir_j.
\tag{5.11}
\]

Let

\[
m_{ij}=\langle C_i,C_j\rangle,
\qquad p_0=a_1-a_2s,
\qquad \theta=\beta-2a_1s+a_2s^2.
\]

Substituting (5.8) into (5.11) and using

\[
B_iB_j^T=m_{ij}+E_iC_j^T+B_iE_j^T
\]

gives

\[
\begin{aligned}
m_{ij}-\theta={}&
p_0(\epsilon_i+\epsilon_j)
+a_2\epsilon_i\epsilon_j
-B_iH(B)_{*,j}\\
&-E_iC_j^T-B_iE_j^T.
\end{aligned}
\tag{5.12}
\]

Consequently,

\[
|m_{ij}-\theta|
\le2(1+p_0)\delta+a_2\delta^2+\sqrt{nU}.
\tag{5.13}
\]

Indeed, the last three terms in (5.12) are bounded by \(\sqrt n\,h_j\),
\(\sigma_i\), and \(\sigma_j\), respectively; moreover, \(h_j^2\le U\)
and \(|\epsilon_i|,|\epsilon_j|\le\delta\).

Direct simplification gives

\[
\theta-\rho=
\begin{cases}
-\dfrac{n}{4(n-1)},&n\equiv0\pmod4,\\[6pt]
\dfrac{n-2}{4(n-1)},&n\equiv2\pmod4.
\end{cases}
\tag{5.14}
\]

Define the rational number

\[
R=1-|\theta-\rho|-2(1+p_0)\delta-a_2\delta^2.
\]

For every parameter pair \((n,t)\) in (5.1), substitution gives

\[
R>0,
\qquad
R^2>nU.
\tag{5.15}
\]

The strict rational lower bounds, after clearing denominators, are as
follows.

\begin{longtable}[]{@{}lrr@{}}
\toprule\noalign{}
Orders & Lower bound for \(R\) & Lower bound for \(R^2-nU\) \\
\midrule\noalign{}
\endhead
\bottomrule\noalign{}
\endlastfoot
\(4\) & \(1/2\) & \(3/50\) \\
\(8\) & \(1/2\) & \(1/25\) \\
\(10\) & \(3/5\) & \(1/100\) \\
\(12\) & \(3/5\) & \(1/25\) \\
\(14\) & \(3/5\) & \(1/50\) \\
\(16,20,\ldots,64\) & \(13/20\) & \(1/20\) \\
\(18,22,\ldots,62\) & \(7/10\) & \(3/40\) \\
\end{longtable}

Since \(R>0\), (5.15) gives \(\sqrt{nU}<R\). Combining this with (5.13),
(5.14), and the definition of \(R\), we obtain

\[
|m_{ij}-\rho|<1.
\]

Since \(m_{ij}-\rho\) is an integer,

\[
\boxed{m_{ij}=\rho}
\qquad(i\ne j).
\]

This proves the lemma. \(\square\)

\hypertarget{exact-finite-projection}{%
\subsection{Exact finite projection}\label{exact-finite-projection}}

We now separate the fixed term \(-s\) from the error \(\epsilon_i\) in
(5.8). Define

\[
\kappa_n=(-1)^{n/2+1}
\left(\beta-2a_1s+a_2s^2-\rho\right).
\]

Using (5.14),

\[
\boxed{
\kappa_n=\frac14+\frac{(-1)^{n/2}}{4(n-1)}.
}
\]

Thus the fixed off-diagonal contribution has magnitude \(\kappa_n\),
while the corresponding fixed term on the distinguished diagonal
cancels.

Let the selected binary rows form \(\widehat C\in\{0,1\}^{m\times n}\),
with the distinguished row indexed by \(\ell\), and put

\[
y=(-1)^{n/2+1}H(B)_{*,\ell},
\qquad
z=\widehat C y.
\]

Lemmas 5.1 and 5.2 imply

\[
\widehat C\widehat C^T=qI_m+\rho J_m.
\tag{5.16}
\]

Since (5.16) is positive definite, the least-norm solution of
\(\widehat Cx=z\) gives

\[
\|y\|_2^2\ge
z^T(qI_m+\rho J_m)^{-1}z,
\tag{5.17}
\]

where

\[
(qI_m+\rho J_m)^{-1}
=q^{-1}I_m-
\frac{\rho}{q(q+m\rho)}J_m.
\tag{5.18}
\]

It remains to verify a finite family of box-projection estimates,
indexed by one-dimensional intervals.

Put

\[
\overline U_\ell=\frac{n(n-2)}{4(n-1)^2}\frac{\tau_n}{n},
\qquad
\Delta(u)=\frac{1-\sqrt{1-4u}}2.
\]

Equations (5.5)--(5.6), applied with \(d=n\), give
\(\overline U_\ell<1/8\), so \(\Delta\) is well defined on
\([0,\overline U_\ell]\).

Recall that every non-distinguished selected row satisfies
\(\sigma_i\le\delta\), with \(\delta\) defined in (5.10). Divide
\([0,\overline U_\ell]\) into \(N=4000\) intervals

\[
I_k=[u_0,u_1]
=\left[\frac{k\overline U_\ell}{N},
\frac{(k+1)\overline U_\ell}{N}\right],
\qquad 0\le k<N.
\tag{5.19}
\]

The mesh size \(N=4000\) was chosen so that every interval has a
positive exact rational margin while keeping the certificate modest; no
optimality is intended.

On such an interval set

\[
\Delta_\ell=\Delta(u_1),
\qquad
\overline h^{\,2}
=\overline U_\ell-u_0-\gamma\max\{0,s-\Delta_\ell\}^2.
\tag{5.20}
\]

Lemma 5.1 applied with \(d=n\) shows that the distinguished row has
weight \(n/2\), and (5.7) gives \(\sigma_\ell<1/2\). If the right-hand
side of (5.20) is negative, the interval contains no feasible
distinguished row. Otherwise, for \(u_\ell\in I_k\), the inequality
\(u_\ell\ge\sigma_\ell(1-\sigma_\ell)\) gives
\(\sigma_\ell\le\Delta(u_1)\). Subtracting the lower bounds
\(u_\ell\ge u_0\) and \(|r_\ell|\ge\max\{0,s-\Delta_\ell\}\) from (5.4)
then gives

\[
\sigma_\ell\le\Delta_\ell,
\qquad
h_\ell^2\le\overline h^{\,2}.
\tag{5.21}
\]

Write \(\overline h=\sqrt{\overline h^{\,2}}\). For \(n\equiv0\pmod4\),
define

\[
L=\kappa_n-
\left[(1+p_0)(\delta+\Delta_\ell)
+a_2\delta\Delta_\ell+\delta\overline h\right].
\tag{5.22}
\]

For \(n\equiv2\pmod4\), retaining the signs of the rounding masses gives
the sharper bound

\[
L=\kappa_n+
\min_{\substack{e\in\{-\delta,\delta\}\\
f\in\{-\Delta_\ell,\Delta_\ell\}}}
\left[
p_0(e+f)+a_2ef
-\frac{\delta+e}{2}
-\frac{\Delta_\ell+f}{2}
-\delta\overline h
\right].
\tag{5.23}
\]

To justify the sign-sensitive terms, decompose the error in any row into
its positive and negative parts. Since \(\epsilon=\sum_jE_j\) and
\(\sigma=\sum_j|E_j|\),

\[
\sum_{E_j>0}E_j=\frac{\sigma+\epsilon}{2}.
\]

Consequently,

\[
E_iC_\ell^T\le\frac{\delta+\epsilon_i}{2},
\qquad
B_iE_\ell^T\le\frac{\Delta_\ell+\epsilon_\ell}{2}.
\]

Rearranging (5.12), using \(m_{i\ell}=\rho\), and applying these two
bounds gives (5.22)--(5.23). The expression in (5.23) is affine in
either variable with the other fixed, so its minimum over
\([-\delta,\delta]\times[-\Delta_\ell,\Delta_\ell]\) occurs at a corner.
Thus

\[
z_i\ge L\qquad(i\ne\ell).
\tag{5.24}
\]

On the diagonal, the fixed-weight identity gives

\[
\|B_\ell\|_2^2
=\sum_j B_{\ell j}^2
=\sum_j B_{\ell j}-u_\ell
=\frac n2+\epsilon_\ell-u_\ell.
\]

The \((\ell,\ell)\) entry of (2.5) is therefore

\[
B_\ell H(B)_{*,\ell}
=\frac{n^2}{4(n-1)}+\beta
+2a_1r_\ell+a_2r_\ell^2
-\left(\frac n2+\epsilon_\ell-u_\ell\right).
\]

After substituting \(r_\ell=-s+\epsilon_\ell\), the constant term
vanishes because

\[
\frac{n^2}{4(n-1)}+\beta-2a_1s+a_2s^2-\frac n2=0,
\]

while the linear coefficient is \(2a_1-2a_2s-1=2p_0-1\). Hence

\[
B_\ell H(B)_{*,\ell}
=u_\ell+(2p_0-1)\epsilon_\ell+a_2\epsilon_\ell^2.
\]

Since

\[
C_\ell H(B)_{*,\ell}
=B_\ell H(B)_{*,\ell}-E_\ell H(B)_{*,\ell},
\]

Cauchy--Schwarz and (5.21) yield

\[
|z_\ell|\le
R_k:=u_1+|2p_0-1|\Delta_\ell
+a_2\Delta_\ell^2+\Delta_\ell\overline h.
\tag{5.25}
\]

The inverse in (5.18) shows that, over the region (5.24)--(5.25), the
quadratic form is minimized at

\[
z_i=L\ (i\ne\ell),
\qquad z_\ell=R_k,
\]

provided

\[
L(q+\rho)>\rho R_k,
\qquad
R_k(q+\rho t)<\rho tL.
\tag{5.26}
\]

Indeed, these are precisely the two inward derivative inequalities for
that corner. At the corner the value is

\[
Q_{n,k}=
\frac{tL^2+R_k^2}{q}
-\frac{\rho}{q(q+(t+1)\rho)}(tL+R_k)^2.
\]

\textbf{Lemma 5.3.} For every parameter pair in (5.1) and every feasible
interval (5.19), the inequalities (5.26) hold and

\[
Q_{n,k}>\overline h^{\,2}.
\tag{5.27}
\]

Lemma 5.3 is verified in the accompanying Lean development. The
calculation uses rational and integer arithmetic, with directed rational
enclosures for the square roots justified by exact squaring. It covers
all thirty parameter pairs and the \(4000\) intervals in (5.19) for each
pair; when \(n\equiv2\pmod4\), it also compares the four corners in
(5.23). An independent implementation in Python gives the same positive
margins. The least certified difference \(Q_{n,k}-\overline h^{\,2}\)
occurs for \((n,t)=(10,5)\). Writing this least value as \(\mu\), its
exact fraction is recorded in the certificate output, and

\[
\mu>
\frac{1003}{2{,}000{,}000}>0.
\]

Section 7 explains how this calculation is connected to the formal
proof.

Since \(\|y\|_2^2=h_\ell^2\le\overline h^{\,2}\), equations (5.17) and
(5.24)--(5.27) give

\[
h_\ell^2=\|y\|_2^2
\ge z^T(qI_m+\rho J_m)^{-1}z
\ge Q_{n,k}
>\overline h^{\,2},
\]

a contradiction.

For each even order

\[
4\le n\le64,
\qquad
n\ne6.
\]

A matrix satisfying (3.1) at any such order would produce the
distinguished row in (5.4), and the preceding argument would give this
contradiction. Hence (3.1) is impossible throughout the stated range.

\hypertarget{the-six-dimensional-case-n6}{%
\section{\texorpdfstring{The six-dimensional case:
\(n=6\)}{The six-dimensional case: n=6}}\label{the-six-dimensional-case-n6}}

At \(n=6\), the one-column inequalities from Section 5 are compatible
with the defect budget. The argument must instead use several columns of
\(H(B)\) at once. We begin with the full identity (2.5) and the global
budget.

\hypertarget{the-defect-budget-and-five-rounded-triples}{%
\subsection{The defect budget and five rounded
triples}\label{the-defect-budget-and-five-rounded-triples}}

At \(n=6\),

\[
r=B\mathbf1-\frac{25}{8}\mathbf1,
\qquad
\gamma=\frac4{25},
\]

and

\[
\mathcal U_i=u_i+h_i^2+\frac4{25}r_i^2.
\]

Dividing (2.6) by \(K=25/6\) gives

\[
\sum_{i=1}^6\mathcal U_i
\le
\frac6{25}\left(\frac1{16}+\frac{33}{49}\right)
=\frac{1731}{9800}
=:\mathcal U.
\tag{6.1}
\]

Separate the non-\(H\) part by setting

\[
P=\sum_{i=1}^6
\left(u_i+\frac4{25}r_i^2\right).
\]

With this notation, (6.1) becomes

\[
\boxed{
P+\|H(B)\|_F^2\le\mathcal U.
}
\tag{6.2}
\]

For every row, (3.3) becomes

\[
r_i=(w_i-3)+\epsilon_i-\frac18.
\]

A row with \(w_i\ne3\) contributes at least

\[
\mathcal U_i\ge\frac{49}{400}.
\tag{6.3}
\]

Indeed, if \(\mathcal U_i<49/400\), then

\[
\sigma_i\le2u_i<\frac{49}{200}<\frac12,
\]

and

\[
u_i\ge\sigma_i(1-\sigma_i),
\qquad
|r_i|\ge\frac78-\sigma_i.
\]

Consequently,

\[
\mathcal U_i\ge
\sigma_i(1-\sigma_i)
+\frac4{25}\left(\frac78-\sigma_i\right)^2.
\]

The right-hand side is concave on \([0,1/2]\), and both endpoint values
are at least \(49/400\), a contradiction. This proves (6.3).

Two rows of the wrong weight would therefore contribute at least

\[
2\cdot\frac{49}{400}
=\frac{49}{200}
>\frac{1731}{9800}.
\]

Therefore at least five rows of \(C\) have weight three. Fix five of
them and, for the remainder of this section, let \(i,j\) range over
their indices. Then

\[
r_i=-\frac18+\epsilon_i.
\tag{6.4}
\]

\hypertarget{exceptional-pair-measurements}{%
\subsection{Exceptional-pair
measurements}\label{exceptional-pair-measurements}}

At order six, (2.5) is

\[
B(B^T+H(B))
=\frac95I+\frac{125}{96}J
+\frac5{12}(r\mathbf1^T+\mathbf1r^T)
+\frac2{15}rr^T.
\tag{6.5}
\]

Let

\[
m_{ij}=\langle C_i,C_j\rangle,
\qquad
T_{ij}=(BH(B))_{ij}=B_iH(B)_{*,j}.
\]

The right-hand side of (6.5) is symmetric, as is \(BB^T\). Consequently
\(BH(B)=B(B^T+H(B))-BB^T\) is symmetric, so \(T_{ij}=T_{ji}\). We retain
the directed notation because the same scalar is measured against two
different columns of \(H(B)\).

For \(i\ne j\), substitute (6.4) into (6.5) and expand \(B=C+E\) to
obtain

\[
\begin{aligned}
T_{ij}={}&
\frac65-m_{ij}
+\frac25(\epsilon_i+\epsilon_j)
+\frac2{15}\epsilon_i\epsilon_j\\
&-E_iC_j^T-C_iE_j^T-E_iE_j^T.
\end{aligned}
\tag{6.6}
\]

Call \(\{i,j\}\) an \textbf{exceptional pair} if \(m_{ij}\ne1\). Since
two triples on six points have intersection \(0,1,2\), or \(3\), every
exceptional pair satisfies

\[
\left|\frac65-m_{ij}\right|\ge\frac45.
\]

Because \(\eta_{ik}\le1/2\),

\[
\|E_i\|_2^2\le u_i.
\]

Suppose \(P\le P_0<1/4\), and choose \(\alpha<1/2\) such that

\[
P_0<\alpha(1-\alpha).
\tag{6.7}
\]

Then every \(\eta_{ik}<\alpha\). Indeed, if some \(\eta_{ik}\ge\alpha\),
the monotonicity of \(x(1-x)\) on \([0,1/2]\) would give

\[
P\ge u_i\ge \eta_{ik}(1-\eta_{ik})
\ge\alpha(1-\alpha)>P_0,
\]

contrary to \(P\le P_0\). It follows also that

\[
\sum_i\sigma_i<\alpha.
\tag{6.8}
\]

Indeed,

\[
\sum_i u_i
=\sum_{i,k}\eta_{ik}(1-\eta_{ik})
\ge(1-\alpha)\sum_i\sigma_i,
\]

whereas

\[
\sum_i u_i\le P\le P_0
<\alpha(1-\alpha).
\]

Since \(1-\alpha>0\), this separate strict inequality gives (6.8).

Using \(|\epsilon_i|\le\sigma_i\), (6.6), (6.8), the
arithmetic--geometric mean inequality, and \(\sum_i u_i\le P_0\), every
selected pair satisfies

\[
\left|
T_{ij}-\left(\frac65-m_{ij}\right)
\right|<\omega(\alpha,P_0),
\tag{6.9}
\]

where

\[
\boxed{
\omega(\alpha,P_0)
=\frac75\alpha
+\frac1{30}\alpha^2
+\frac12P_0.
}
\tag{6.10}
\]

The three terms in (6.10) respectively bound

\[
\frac25(\epsilon_i+\epsilon_j)-E_iC_j^T-C_iE_j^T,
\qquad
\frac2{15}\epsilon_i\epsilon_j,
\qquad
E_iE_j^T.
\]

Moreover, every selected row satisfies

\[
\|B_i\|_2^2=3+\epsilon_i-u_i<3+\alpha.
\tag{6.11}
\]

\hypertarget{the-energy-contradiction}{%
\subsection{The energy contradiction}\label{the-energy-contradiction}}

Five triples on a six-point set cannot be pairwise intersecting in
exactly one point. If \(d_k\) is the number of triples containing point
\(k\), pairwise intersection one would imply

\[
\sum_{k=1}^6d_k=15,
\qquad
\sum_{k=1}^6\binom{d_k}{2}=\binom52=10,
\]

and hence

\[
\sum_{k=1}^6d_k^2=35.
\]

But Cauchy--Schwarz gives

\[
\sum_{k=1}^6d_k^2
\ge\frac{15^2}{6}=\frac{75}{2}>35.
\]

Thus there is at least one exceptional pair.

If there is exactly one exceptional pair, its two triples must be equal.
Indeed, intersection zero is impossible because the two triples are
complementary, so no third triple can intersect both in one point. If
the intersection is two, write the triples as \(\{a,b,x\}\) and
\(\{a,b,y\}\). A triple meeting both in one point is one of

\[
\{a,z,w\},
\quad
\{b,z,w\},
\quad
\{x,y,z\},
\quad
\{x,y,w\}.
\]

The first two triples intersect in \(\{z,w\}\), and the last two
intersect in \(\{x,y\}\). Any choice of three therefore contains one of
these pairs and cannot have all pairwise intersections equal to one.
Hence a unique exceptional pair must have \(m_{ij}=3\), and its
unperturbed signal has magnitude \(9/5\).

Suppose instead that there are at least two exceptional pairs. Choose
two. If they are disjoint, their two orientations use four distinct
columns of \(H(B)\). If they meet, say \(\{i,j\}\) and \(\{i,k\}\),
columns \(j\) and \(k\) each give one scalar constraint, while column
\(i\) gives two jointly.

More precisely, suppose that every chosen measurement has magnitude at
least \(a\) and every selected row has squared norm at most \(V\).
Cauchy--Schwarz and the two-row operator-norm bound give

\[
\boxed{
\|H(B)\|_F^2\ge\frac{3a^2}{V}.
}
\tag{6.12}
\]

In the disjoint case, Cauchy--Schwarz gives a contribution at least
\(a^2/V\) from each of four distinct columns, which is more than (6.12)
requires. For the shared column in the second case, let \(M\) be the
matrix consisting of the two relevant rows of \(B\). Then

\[
\|MH(B)_{*,i}\|_2^2
\le\|M\|_2^2\|H(B)_{*,i}\|_2^2
\le2V\|H(B)_{*,i}\|_2^2,
\]

whereas the left-hand side is at least \(2a^2\). The shared column
therefore contributes at least \(a^2/V\); the two remaining columns
contribute another \(a^2/V\) each. This proves (6.12).

There are two ranges for \(P\).

The breakpoint \(1/25\) and the two values of \(\alpha\) below were
chosen together: they keep the rounding estimate valid in each range and
leave a positive exact energy margin on both sides of the breakpoint.

\hypertarget{case-a-pge125}{%
\subsubsection{\texorpdfstring{Case A:
\(P\ge1/25\)}{Case A: P\textbackslash ge1/25}}\label{case-a-pge125}}

Take

\[
\alpha=\frac{23}{100},
\qquad
P_0=\mathcal U=\frac{1731}{9800}.
\]

Condition (6.7) holds, since

\[
\frac{23}{100}\frac{77}{100}-\mathcal U
=\frac{229}{490000}>0.
\]

Moreover,

\[
\omega\left(\frac{23}{100},\mathcal U\right)
=\frac{6057571}{14700000}
<\frac{207}{500}.
\tag{6.13}
\]

Thus every exceptional-pair measurement satisfies

\[
|T_{ij}|>\frac{193}{500},
\]

and (6.11) gives

\[
\|B_i\|_2^2<\frac{323}{100}.
\]

If there are at least two exceptional pairs, (6.12) and \(P\ge1/25\)
imply

\[
\begin{aligned}
P+\|H(B)\|_F^2
&>\frac1{25}
+3\frac{(193/500)^2}{323/100}\\
&=\frac{144047}{807500}
>\frac{1731}{9800},
\end{aligned}
\]

where the last margin is

\[
\frac{138781}{79135000}>0.
\]

This contradicts (6.2).

If there is exactly one exceptional pair, it is a duplicate pair. By
(6.9) and (6.13), its two directed measurements have magnitude greater
than

\[
\frac95-\frac{207}{500}=\frac{693}{500}.
\]

Consequently,

\[
\|H(B)\|_F^2
>2\frac{(693/500)^2}{323/100}
=\frac{480249}{403750}
>\mathcal U,
\]

again contradicting (6.2).

\hypertarget{case-b-p125}{%
\subsubsection{\texorpdfstring{Case B:
\(P<1/25\)}{Case B: P\textless1/25}}\label{case-b-p125}}

Take

\[
\alpha=\frac{21}{500},
\qquad
P_0=\frac1{25}.
\]

Again (6.7) holds, because

\[
\frac{21}{500}\frac{479}{500}-\frac1{25}
=\frac{59}{250000}>0,
\]

and

\[
\omega\left(\frac{21}{500},\frac1{25}\right)
=\frac{197147}{2500000}
<\frac2{25}.
\]

Every exceptional-pair measurement therefore has magnitude greater than
\(18/25\), and every selected row has squared norm less than
\(1521/500\). If there are at least two exceptional pairs, (6.12) gives

\[
\|H(B)\|_F^2
>3\frac{(18/25)^2}{1521/500}
=\frac{432}{845}
>\mathcal U,
\]

contradicting (6.2). If there is exactly one exceptional pair, it is a
duplicate pair, and its two directed measurements have magnitude greater
than

\[
\frac95-\frac2{25}=\frac{43}{25}.
\]

Therefore

\[
\|H(B)\|_F^2
>2\frac{(43/25)^2}{1521/500}
=\frac{14792}{7605}
>\mathcal U,
\]

again contradicting (6.2).

Thus every possible order-six configuration contradicts (6.2), and hence

\[
\|B^{-1}\|_F>\frac{12}{7}.
\]

The use of several columns is essential here. A single distinguished
column sees only the weak \(1/5\) signal, whereas two exceptional pairs
distribute enough energy among at least three columns of \(H(B)\) to
exceed the global budget once the row defects are included.

\hypertarget{formal-verification}{%
\section{Formal verification}\label{formal-verification}}

The even-dimensional proof has been formalized in Lean 4.19.0 with
mathlib 4.19.0. The development begins with an arbitrary nonsingular
nonnegative matrix and includes the normalization, the global defect
estimate, the rounding arguments in Sections 4 and 5, the order-six
argument, and the direct calculation for order two.

At the level of the paper, the imported content of Frankel--Urschel
Lemma 2.1 consists of the estimates (2.3) and (2.4); identity (2.5) was
derived above from the explicit definition (2.2). The Lean interface is
deliberately more abstract: it retains witnesses satisfying all three
conclusions of the published lemma, including (2.5), and uses no other
property of \(H(B)\). Its hypotheses, constants, indices, and transposes
were checked against the published statement. Thus the formalization is
conservatively conditional on {[}3, Lemma 2.1{]}, whose proof has not
itself been formalized.

The finite calculation in Section 5 uses rational and integer
arithmetic, without floating-point approximations. Its Lean
implementation includes a kernel-checked soundness theorem showing that
a successful check implies the inequalities used in the matrix argument.
The Boolean equality asserting that the checker succeeds is discharged
by \texttt{native\_decide}; this evaluation step uses the standard
bridge \texttt{Lean.ofReduceBool} and therefore trusts Lean's native
evaluator and compiler in addition to the kernel. Source code, exact
margins, and instructions for reproducing the build accompany the paper.
The separate Python implementation is an independent regression check
and is not imported into Lean.

\hypertarget{completion-of-the-proof}{%
\section{Completion of the proof}\label{completion-of-the-proof}}

It remains only to combine the even-dimensional argument with the
elementary small orders and Cheng's theorem.

For \(n=1\), write \(A=(a)\) with \(a>0\). Then

\[
\|A^{-1}\|_F=a^{-1}=\|A\|_{\max}^{-1},
\]

which is Theorem 1.1 with equality; the normalized matrix \((1)\) is the
unique \(S\)-matrix of order one.

For \(n=2\), normalize

\[
B=\|A\|_{\max}^{-1}A=
\begin{pmatrix}a&b\\c&d\end{pmatrix},
\qquad 0\le a,b,c,d\le1.
\]

Since \(B\) is nonsingular,

\[
\|B^{-1}\|_F^2
=\frac{a^2+b^2+c^2+d^2}{(ad-bc)^2}.
\]

If \(ad\ge bc\), nonsingularity gives \(0<ad-bc\le ad\), so in
particular \(ad>0\), and hence

\[
\|B^{-1}\|_F^2
\ge\frac{a^2+d^2}{(ad)^2}
\ge\frac2{ad}\ge2.
\]

If \(bc>ad\), the same argument with \(b,c\) gives the same conclusion.
Therefore

\[
\|A^{-1}\|_F
=\|A\|_{\max}^{-1}\|B^{-1}\|_F
\ge\sqrt2\,\|A\|_{\max}^{-1}
>\frac43\,\|A\|_{\max}^{-1}.
\]

For every odd \(n\ge3\), Theorem 1.1, including its equality
characterization, is Cheng's theorem {[}1{]}. For even \(n\ge4\), the
contradiction in Section 4 covers \(n\ge66\), the calculation in Section
5 covers \(4\le n\le64\) with \(n\ne6\), and Section 6 covers \(n=6\).
Hence the conjectured inequality holds in every dimension. Equality is
possible in the odd-dimensional case exactly for positive multiples of
\(S\)-matrices, while the inequality is strict in every even order. This
proves Theorem 1.1. \(\square\)

The proof above recovers the S-matrix theorem by a different mechanism
from Zhang's centered spectral argument. Here the contradiction passes
through nearly integral rows, exact intersection numbers, and column energy;
the finite and order-six branches are retained as part of the independent
proof. The formal verification provides an additional check of this route,
with the external mathematical input and native-evaluation boundary stated
in Section 7.

\hypertarget{code-and-data-availability}{%
\section*{Code and data availability}\label{code-and-data-availability}}
\addcontentsline{toc}{section}{Code and data availability}

The Lean source, the exact-rational Python certificate, all recorded
margins, and complete reproduction instructions are supplied as supplementary files with this manuscript and are also available
in the ancillary archive of arXiv:2608.29750. The distribution fixes Lean and
mathlib at version 4.19.0 through \texttt{lean-toolchain} and
\texttt{lakefile.toml}; the supplied preflight script rebuilds the Lean
project, reruns the Python certificate, and rejects proof placeholders
or project-defined axioms.

\hypertarget{computational-and-editorial-assistance}{%
\section*{Computational and editorial
assistance}\label{computational-and-editorial-assistance}}
\addcontentsline{toc}{section}{Computational and editorial assistance}

Computer algebra and exact-arithmetic programs were used to explore the
finite inequalities. The final argument was implemented in Lean with a
proved checker-to-matrix soundness theorem. Generative AI tools assisted
with exploratory programming, formalization, and editorial review. The
author checked all mathematical claims, constants, citations, and source
files; the scope and trust boundary of the formal verification are
stated in Section 7.

\hypertarget{references}{%
\section*{References}\label{references}}
\addcontentsline{toc}{section}{References}

\leavevmode\par

\begin{enumerate}
\def\labelenumi{\arabic{enumi}.}
\tightlist
\item
  Ching-Shui Cheng, ``An Application of the Kiefer--Wolfowitz
  Equivalence Theorem to a Problem in Hadamard Transform Optics,''
  \emph{The Annals of Statistics} \textbf{15} (1987), no. 4, 1593--1603.
  https://doi.org/10.1214/aos/1176350612.
\item
  Roman Drnovšek, ``On the S-matrix Conjecture,'' \emph{Linear Algebra
  and its Applications} \textbf{439} (2013), no. 11, 3555--3560;
  arXiv:1306.6786.
\item
  Elsa Frankel and John Urschel, ``On the Frobenius Norm of the Inverse
  of a Non-Negative Matrix,'' \emph{Linear Algebra and its Applications}
  \textbf{708} (2025), 193--203.
  https://doi.org/10.1016/j.laa.2024.11.030; arXiv:2409.04354.
\item
  N. J. A. Sloane and Martin Harwit, ``Masks for Hadamard Transform
  Optics, and Weighing Designs,'' \emph{Applied Optics} \textbf{15}
  (1976), 107--114.
\item
  Martin Harwit and N. J. A. Sloane, \emph{Hadamard Transform Optics},
  Academic Press, New York, 1979.
\item
  Teng Zhang, ``The A-optimal Design Problem and the S-matrix Conjecture,''
  arXiv:2609.20306v1 (2026).
  \url{https://arxiv.org/abs/2609.20306}.
\end{enumerate}

\end{document}